\documentclass{article}
\usepackage{amsmath}
\usepackage{amssymb}
\usepackage[utf8]{inputenc}
\usepackage[english]{babel}
\usepackage{graphicx}
\usepackage{theorem}
\usepackage{subcaption}
\usepackage{algorithm}
\usepackage{algorithmic}
\usepackage{color}
\usepackage{fullpage}
\title{Swimming at Low Reynolds number}
\newcommand{\FiguresPath}{.}
\author{Luca Berti, Laetitia Giraldi, Christophe Prud'homme}
\begin{document}

	%
	\maketitle
	
	%
	%
	\begin{paragraph}{Abstract.}
		We address the swimming problem at low Reynolds number. This regime, which is typically used for micro-swimmers, is described by Stokes equations. We couple a PDE solver of Stokes equations, derived from the Feel++ finite elements library, to a quaternion-based rigid-body solver. We validate our numerical results both on a 2D exact solution and on an exact solution for a rotating rigid body respectively. Finally, we apply them to simulate the motion of a one-hinged swimmer, which obeys to the scallop theorem.
	\end{paragraph}
	\begin{paragraph}{Resumé.}
		Nous considérons le problème de la nage à bas nombre de Reynolds. Ce régime, typique des micro-organismes nageurs, est décrit par les équations de Stokes. Nous couplons un solveur d'EDP des équations de Stokes, construit à l'aide de la librairie aux éléments finis Feel++, avec un solveur de corps rigide basé sur les quaternions. Nous validons les deux solveurs à l'aide d'une solution exacte en 2D pour le fluide et d'une solution exacte pour un corps rigide qui tourne. Nous les appliquons pour simuler la nage d'un micro-organisme qui obéit au théorème de la coquille de Saint-Jacques.
	\end{paragraph}
	\maketitle
	\section*{Introduction}
	One of the first works in the field of microswimming was proposed by E. Purcell, who presented in 1977 his observations regarding micro-swimmers' swimming strategies. In \cite{Purcell1977} he noticed that, at the microscale, Reynolds number is small, which means that viscous dissipative forces overcome inertial ones. This observation has two consequences: firstly, Navier-Stokes fluid equations, at the swimmer's scale, can be approximated very well by Stokes equations. These have the characteristics of instantaneously diffusing any change in momentum, and they also guarantee the adherence of the fluid to the boundary of the swimmer at every instant; secondly, certain types of deformations, that usually at human scale produce a net displacement, do not give any motion at the micro-scale. The well-known Scallop theorem \cite{Chambrion2010,Purcell1977} specifies some of them.
	
	Since Stokes equations are time-invariant, the important factors, in terms of swimming, are just the body configurations and not the speed of deformation. 
	
	A typical value of Reynolds number $Re = \frac{UL}{\nu}$ at the microscale is $10^{-5}$. In this case the length-scale of microscopic swimmers is taken to be $L \approx 10^{-6}$ $m$, their mean  velocity $U\approx 30$ $\mu m/s$ and water kinematic viscosity $\nu \approx 10^{-6}$ $m^2/s$.
	
	Mathematical analysis of microswimming has mainly developed in the field of control theory, due to its applications to micro-robotics. In fact, it is a matter of great importance to understand how the propelling mechanisms of a robot influence its motion and its controllability. As microrobots might be used as vessels for targeted drug administration or as microsurgery devices, the necessity of control theory becomes clear. In \cite[chapter 16]{Rosen2011} some applications are mentioned: they can range from medical imaging of regions that cannot be reached through endoscopy, to targeted drug delivery or micro-surgery. Different propulsion techniques can be devised: in \cite{Ahmed2017} a combined acoustic and magnetic control was proposed, while in \cite{Oulmas2017} a magnetic control alone was studied.
	
	Microswimming (under the approximation given by the Stokes regime) has already been analyzed from the ODE point of view using Resistive Force Theory \cite{GRAY802}. This theory allows the approximation of fluid effects via appropriate coefficients. Such coefficients, inserted in the ordinary differential equations that arise from the Newton equations of motion, approximate horizontal and vertical drag forces for a variety of swimmers made of simple shapes, like segments or circles. The approximation deriving from Resistive Force Theory makes it possible to use ODE control techniques in this framework \cite{Alouges2013,Alouges2015,Pomet2017,Addendum2018,Martinon2015}. 
	
	Swimming simulations have been addressed also in the PDE setting by several authors studying fluid-structure interaction. In this case, fluid forces could be computed for swimmers with different shapes. Several techniques can be applied, like level set methods \cite{IolloBergmann2011,IolloBergmann2016}, fictitious domain methods \cite{Aguillon2013} or immersed boundary methods \cite{FAUCI1995679}. Other schemes, like boundary element methods, can be applied when the governing equations have an associated Green function, like Stokes problem does (see \cite{PimponiChinappi2018}).

	Up to now, literature has mainly considered microswimmers that can be represented by finite dimensional control systems, and very often the effects of fluid forces on the body and of its elasticity are modeled in a simplified manner, in order to still use tools from control theory of ODEs. Swimmers of different shapes have been investigated: among them there are Purcell 3-link microswimmer, the N-link microswimmer \cite{Alouges2013}, sphere-based microswimmers \cite{DesimoneMontino2015} and some other models that add elasticity effects through torsional springs or control the swimmer through an externally imposed magnetic field. An important challenge, coming from microrobot control and linked to the applications we have mentioned earlier, is controlling the object close to walls or obstacles. Some results concerning controllability in the half space have been obtained by L. Giraldi \cite{AlougesGiraldi2012,VaretGiraldi2013}.
	
	Other studies about control of microswimmers have been presented in \cite{lohea2012}, where a deformable amoeba was studied, and in subsequent papers by Lohéac and Munnier \cite{loheac_munnier}. Other studies focused on the geometry arising from microswimming, for example \cite{Chambrion2010, ShapereWilczek1989}. In \cite{Chambrion2010} the scallop theorem is extended using Riemannian geometry and in \cite{giraldi2017} optimal control of a finite-dimensional swimmer moving on a straight line is studied.

	The aim of this paper is to compute the displacement of a microscopic scallop in order to illustrate the Scallop theorem of Purcell, using the finite element library Feel++.

	The paper is organised as follows: in Section \ref{Section:Model} we will present the mathematical model and the governing equations of swimmer's displacement. Section \ref{Section:Method} introduces the associated discrete problem. Section \ref{Section:Implementation} gives some details about the implementation of it. Section 	\ref{Section:Experiments} presents our numerical results on the scallop swimmer. A paragraph with remarks and perspectives concludes the paper.
	
	\section{The model}
	\label{Section:Model}
	In this section, the equations that are considered along the paper are introduced. Section \ref{Subsection:Stokes} is devoted to the description of the fluid model. Section \ref{Subsection:RigidBody} describes how the rigid body motion can be recovered given an imposed swimmer body deformation.
	\subsection{The fluid equations}
	\label{Subsection:Stokes}
	In our problem we consider a micro-swimmer immersed in an incompressible fluid which is modeled, in the Eulerian reference frame, by Navier-Stokes equations
	\begin{equation*}
	\begin{aligned}
	&\frac{\partial u}{\partial t} + (u \cdot \nabla)u = \nu \Delta u - \frac{\nabla p}{\rho} + f, \\
	&\nabla \cdot u = 0.
	\end{aligned}
	\end{equation*}
	where $u(t,x):\mathbb{R}\times\mathbb{R}^d \rightarrow \mathbb{R}^d$, $d=2,3$,  is the fluid velocity field, $p(t,x):\mathbb{R}\times\mathbb{R}^d \rightarrow \mathbb{R}$ is the fluid pressure field, $\nu$ is the fluid kinematic viscosity, $\rho$ is the fluid density and $f(t,x):\mathbb{R}\times\mathbb{R}^d \rightarrow \mathbb{R}^d$ is the distribution of external forces per unit mass. The Navier-Stokes system, in the case of micro-swimmers, can be simplified to the Stokes equations by neglecting the inertial contributions of the flow.
	
	 This simplification is justified by the small value of the Reynolds number associated to the flow regime in which we are interested. This can be seen in the adimensionalized version of Navier-Stokes equations (here without external forces, for simplicity)
	\begin{equation*}
	\begin{aligned}
	&\mathrm{Re}\Big(\frac{\partial u'}{\partial t'} + (u' \cdot \nabla)u'\Big) = \Delta u' - \nabla p', \\
	&\nabla \cdot u' = 0,
	\end{aligned}
	\end{equation*} 
	where $u' = \frac{u}{U}$, $p' = \frac{pL}{\mu U}$, $t' = \frac{L}{U}$ and $U$ is the swimmer's characteristic swimming speed, $L$ is its characteristic length and $\nu$ is the fluid kinematic viscosity. In the momentum conservation equations, Reynolds number multiplies the inertial terms, which can therefore be neglected, since in the case of micro-swimmers $\mathrm{Re}\approx10^{-5}$.
	
	As a result, the set of equations we will be using to model the fluid behaviour is the Stokes system
	\begin{equation*}
	\left\{
		\begin{aligned}
		&\mu\Delta u - \nabla p  = 0,\\
		&\nabla \cdot u = 0.
		\end{aligned}
		\right.
	\end{equation*}
	Denoting by $\sigma(t,x) = -p(t,x)\mathbb{I} + \mu(\nabla u(t,x) + \nabla u(t,x)^T)$ the Newtonian fluid stress tensor and considering a time-dependent family of open sets $\Omega^t \subseteq \mathbb{R}^d \, \forall t$, we can therefore write the fluid problem as follows: find $u(t,x),p(t,x)$ such that
	\begin{equation}
		\begin{aligned}
		& \nabla \cdot \sigma = 0 \qquad &\textrm{in $ \mathbb{R}_+\times\Omega^t$},\\
		&\nabla \cdot u = 0 \qquad &\textrm{in $\mathbb{R}_+\times \Omega^t$},\\
		& u(t,x) = \bar{u}(t,x) \qquad &\textrm{on $\mathbb{R}_+\times\partial \Omega^t_D$},
		\\
		& \Big(-p(t,x)\mathbb{I} + \mu \big(\nabla u(t,x) + \nabla u(t,x)^T\big)\Big) \vec{n}(t,x) = \bar{g}(t,x) \qquad &\textrm{on $\mathbb{R}_+\times\partial \Omega^t_N$},
		\end{aligned}
		\label{eq:Stokes}
	\end{equation}
	where $\partial\Omega^t_D$ and $\partial \Omega^t_N$ are the portions of $\partial \Omega^t$ where Dirichlet and, respectively, Neumann boundary conditions are imposed, and $\vec{n}(t,x)$ is the outward unit normal to $\partial \Omega^t$. In particular, one requires $\partial \Omega^t = \partial\Omega^t_D \cup \partial \Omega^t_N$ and $ \partial\Omega^t_D \cap \partial \Omega^t_N = \emptyset$. Conditions on $\bar{u}$ and $\bar{g}$ are imposed so that the problem has a unique solution in $[\mathrm{H}^1(\Omega^t)]^d\times \mathrm{L}^2(\Omega^t)$, notably $\bar{u} \in \textrm{H}^{1/2}(\partial \Omega^t_D)$ and $\bar{g} \in \textrm{H}^{-1/2}(\partial \Omega^t_N)$ provide such result \cite{BrezziFortin1991}.
	
	In the case of our interest, $\Omega^t$ will be $\mathbb{R}^d$ deprived of the region occupied by the swimmer, and the Dirichlet boundary conditions correspond to the velocity field on the swimmer boundary. Such vector field can be decomposed in two parts: one represents the rigid motion and the other comes from the deformation of the swimmer's body $\bar{u}(t,x)=V+\omega \wedge (x-x_0) + v(t,x)$, where $x_0$ is the center of mass of the swimmer. When this latter is known in advance, it is possible to recover the linear velocity $V$ and the angular velocity $\omega$ by solving a linear system as presented in \cite{Munnier2012}. This method exploits the linearity of Stokes equations with respect to pressure and velocity, producing a system of equations which is linear in $(V,\omega)$ and that depends on the stresses produced by a number of ``fundamentally-driven" flows.
	
	\subsection{The rigid body dynamics}
	\label{Subsection:RigidBody}
	Following Newton's second law of dynamics and Euler's equations for angular momentum, we can write the equations describing the dynamics of the swimmer. However, in our case, the microscopic nature of the immersed object and the Stokes approximation of fluid equations lead us to neglect the inertial terms. As a result, 
	the aforementioned equations read respectively (see \cite{Munnier2012})
	\begin{equation*}
	\begin{aligned}
		&\int_{\partial \Omega^t} \sigma(t,x) \vec{n}(t,x) \, \mathrm{d}S(x)= 0, \\
		&\int_{\partial \Omega^t} \sigma(t,x)\vec{n}(t,x)\wedge x \, \mathrm{d}S(x) = 0.
	\end{aligned}
	\label{Equation:Newton}
	\end{equation*}
The equations above are written in the swimmer's reference frame. As proposed and developed in \cite{Munnier2012}, the linearity of the Stokes equation in pressure and velocity allows a decomposition of the two fields as a sum of solutions to fundamentally driven Stokes equations. Each of these problems is characterized by the boundary conditions on the swimmer's surface. If we denote the canonical basis of the swimmer's reference frame by $\{e_i\}_{i=1}^3$ and by $v(t,x)$ the deformation velocity we already mentioned, the fundamental boundary conditions are of the form $\bar{u}(t,x) = e_i$ for $i \in \{1,2,3\}$, $\bar{u}(t,x) = e_i\wedge x$ for $i \in \{1,2,3\}$ and $\bar{u}(t,x) = v(t,x)$. They provide seven different Stokes problems and allow an expansion of the velocity and pressure fields as $u(t,x) = \sum_{i=1}^3V_i u_i(t,x) + \sum_{i=4}^6\omega_iu_i(t,x) + u_7(t,x)$ and $p(t,x) = \sum_{i=1}^3V_i p_i(t,x) + \sum_{i=4}^6\omega_ip_i(t,x) + p_7(t,x)$, where problems 1 to 3 are characterized by $\bar{u}(t,x) = e_i$, problems 4 to 6 are characterized by $\bar{u}(t,x)=e_i\wedge x$ and problem 7 by $\bar{u}(t,x) = v(t,x)$. This decomposition carries onto the fluid stress tensors, since they are linear with respect to pressure and velocity: $\sigma(u,p) = \sum_{i=1}^{3}V_i\sigma(u_i,p_i) + \sum_{i=4}^6\omega_i\sigma(u_i,p_i) + \sigma(u_7,p_7)$. Substituting this decomposition into \eqref{Equation:Newton} and appropriately rearranging the terms brings us to the solution of a linear system in $(V,\omega)$ of $6\times 6$ matrix $M$ and right-hand side $N \in \mathbb{R}^6$ as follows
\begin{equation}
		M_{ij} = 
		\left\{
		\begin{aligned}
		& \int_{\partial \Omega^t} e_i \cdot \sigma(u_j,p_j)\vec{n} \, \mathrm{d}S \qquad &\textrm{for $i\in\{1,2,3\},j\in\{1,\dots,6\}$},\\
		& \int_{\partial \Omega^t} e_{i-3} \cdot \big( \sigma(u_{j},p_{j})\vec{n}\wedge x \big) \, \mathrm{d}S  \qquad &\textrm{for $i\in\{4,5,6\},j \in \{1,\dots,6\}$},
		\end{aligned}
		\right.
		\label{Eq:M}
\end{equation}
\begin{equation}
	N_{i} = 
	\left\{
	\begin{aligned}
	& -\int_{\partial \Omega^t} e_i \cdot \sigma(u_7,p_7)\vec{n} \, \mathrm{d}S  \qquad &\textrm{for $i\in\{1,2,3\}$},\\
	&  -\int_{\partial \Omega^t} e_{i-3} \cdot \big(\sigma(u_{7},p_{7})\vec{n} \wedge x \big)\, \mathrm{d}S \qquad &\textrm{for $i\in\{4,5,6\}$}.
	\end{aligned}
	\right.
	\label{Eq:N}
\end{equation}
	Once the velocities (in the swimmer's reference frame) are recovered by solving the previous system, as done in \cite{Munnier2012}, we can determine the position and orientation of the swimmer in the laboratory reference frame by solving a set of ODEs of the form \cite{Featherstone2007}
	\begin{equation}
		\left\{
		\begin{aligned}
		& \frac{\mathrm{d} q}{\mathrm{d} t}  = \frac{1}{2} [0,R(q)\omega]*q, \\
		& \frac{\mathrm{d} X}{\mathrm{d} t}  = R(q)V,
		\end{aligned}
		\right.
		\label{RigidMotion}
	\end{equation}
	where $q$ is the unit quaternion representing the orientation of the body in the laboratory fixed frame, $R(q)$ is the rotation matrix associated to the unit quaternion $q$, $*$ denotes the quaternion product and $X$ is the position of the swimmer's centroid in the fixed reference frame. Appropriate initial conditions for the ODEs have to be chosen. The derivation of the quaternion equation and definition of the quaternion product are described in section \ref{Section:Method}.
	
	Choosing unit quaternions to represent rotations in $\mathbb{R}^d$ is dictated by the fact that they cover $\mathrm{SO} (3)$, sending the pair $\{q,-q\}$ to a unique $R(q) \in \mathrm{SO}  (3)$, and avoid the classical complications that arise when Euler angles are used \cite{Featherstone2007,Landau1976}.
	\subsection{The self-propulsion constraints}
	In the previous section we addressed the computation of the linear and angular velocity of the swimmer in its reference frame. In order to model the displacement of a swimmer, the deformation velocity $v(t,x)$ has to satisfy two constraints ensuring that the motion is entirely due to internal forces and is not relying on external actions. In particular, these constraints express the fact that in the swimmer's reference frame, linear and angular momenta remain unchanged because no external forces are applied \cite{Chambrion2009}. In the linear momentum case, such invariance writes
	\begin{equation}
		\frac{\mathrm{d}}{\mathrm{d}t}\int_{\Omega^0} x \, dm(x) = 0.
		\label{eq:Chambrion}
	\end{equation}
	Following \cite{Chambrion2009}, \eqref{eq:Chambrion} and its angular momentum analogue can be rewritten as
	\begin{equation*}
		\int_{\partial \Omega^t} v(t,x) \, \mathrm{d}S(x) = 0, \qquad \qquad \int_{\partial \Omega^t} v(t,x) \wedge x \, \mathrm{d}S(x) = 0.
	\end{equation*}	
	\section{Numerical method}
	\label{Section:Method}
	In subsection \ref{Subsection:Numerical_Fluid} we address the discretization of the Stokes system and in subsection \ref{Subsection:NumericalRBD} we will address the derivation and discretization of the ODE system for rigid body motion.
	\subsection{The fluid model}
	\label{Subsection:Numerical_Fluid}
	This part presents the PDE discretization that is used in the rest of the paper to numerically recover the swimmer's motion using finite element method. The indices denoting the time dependence of the fluid domain are abandoned and the focus is placed on the stationary problem to be solved at each time instant. The solution of the Stokes equations was performed via mixed finite elements, where the approximation spaces were chosen in order to guarantee the inf-sup stability of the discrete problem. The variational formulation is the following: find $(u,p) \in [\mathrm{H}^1(\Omega)]^d\cap \{v=\bar{u}|_{\partial \Omega_D}\} \times \mathrm{L}^2(\Omega)$ such that
	\begin{equation*}
	\begin{aligned}
		&\mu\int_{\Omega} (\nabla u+\nabla u ^T) : \nabla \phi \, \mathrm{d}x -\int_{\Omega} p(\nabla \cdot \phi) \, \mathrm{d}x = \int_{\partial \Omega_N} g \cdot \phi \, \mathrm{d}S(x), \\
		&\int_{\Omega} q (\nabla \cdot u) \, \mathrm{d}x= 0,
		\end{aligned}
	\end{equation*}
	for all $(\phi,q) \in [\mathrm{H}^1_0(\Omega)]^d\cap \{v=0|_{\partial \Omega_D}\} \times \mathrm{L}^2(\Omega)$ such that the boundary conditions in \eqref{eq:Stokes} are satisfied.
	
	For the finite element translation, we chose the pair of Taylor-Hood spaces $(X_h^k,M_h^{k-1})$ for velocity and pressure, where
	\[
	\begin{aligned}
	X_h^k &= \{u_h \in [\mathrm{C}^0(\Omega)]^d|\,\,u_h \textrm{ is piecewise polynomial of order $k$ on each element of the triangulation}\}\\
	M_h^{k-1} &= \{p_h \in \mathrm{C}^0(\Omega)|\,\,p_h \textrm{ is piecewise polynomial of order $k-1$ on each element of the triangulation}\}
	\end{aligned}
	\]
	provided that $k\in \mathbb{N},\, k>1$. 
	They are characterized by different approximation orders on $u$ and $p$ in order to guarantee the inf-sup stability of the discrete problem. Namely,
	a continuous, piecewise polynomial approximation of order $k$ is performed on each velocity component and a continuous, piecewise polynomial approximation of order $k-1$ is performed on the pressure.
	
	 In the following, we assume that the computational triangulated domain coincides with the continuous one. $X_h^k$ and $M_h^{k-1}$ are subsets of the functional spaces in which $(u,p)$, solutions of the continuous problem, are chosen. We can write the discrete version of the finite element problem  as follows: find $(u_h,p_h) \in X_h^k\times M_h^{k-1}$ such that 
	\begin{equation*}
		\begin{aligned}
		&\mu\int_{\Omega} (\nabla u_h+\nabla u_h^T) : \nabla \phi_h \, \mathrm{d}x -\int_{\Omega} p_h(\nabla \cdot \phi_h) \, \mathrm{d}x = \int_{\partial \Omega_N} g \cdot \phi_h \, \mathrm{d}S(x), \\
		&\int_{\Omega} q_h (\nabla \cdot u_h) \, \mathrm{d}x= 0,
		\end{aligned}
	\end{equation*}
	for all $(\phi_h,q_h) \in X_h^k\times M_h^{k-1}$. Dirichlet boundary conditions are applied strongly at the algebraic level.
	\subsection{The rigid body dynamics}
	\label{Subsection:NumericalRBD}
	As shown in subsection \ref{Subsection:RigidBody}, from the finite element solution of the fluid problem it is possible to recover the rigid body dynamics of the swimmer. Unit quaternions are a good option to represent the orientation of a rigid object, since they avoid parametrization issues that can affect other representations, like the gimbal lock for Euler angles \cite{Featherstone2007,Landau1976}. 
	
	A quaternion $q$ is defined by a scalar $q_0$ and a 3-dimensional vector $q_v=(q_{v_1},q_{v_2},q_{v_3})$, $q=(q_0,q_v)$. The conjugate of $q$ will be denoted by $\bar{q}=(q_0,-q_v)$ and the quaternion product is defined as $q*r = (q_0r_0-q_v\cdot r_v\, , \,q_0r_v+q_vr_0+q_v\wedge r_v )$. The norm of a quaternion is defined as $||q||=\sqrt{q\bar{q}} = \sqrt{q_0^2+q_{v_1}^2+q_{v_2}^2+q_{v_3}^2}$. Each quaternion of unit norm  represents a reflection in 3D space, and together with its conjugate concurs to describe a rotation in space. The rotation matrix linked to a unit quaternion $q = (q_0,q_v)$ is given by
	\begin{equation*}
	R(q) = 
	\begin{bmatrix}
	q_0^2+q_{v_1}^2-q_{v_2}^2-q_{v_3}^2 & 2(q_{v_1}q_{v_2}-q_0q_{v_3}) & q_{v_1}q_{v_3}+q_0q_{v_2}\\
	2(q_{v_1}q_{v_2}+q_0q_{v_3}) & q_0^2-q_{v_1}^2+q_{v_2}^2-q_{v_3}^2 & 2(q_{v_2}q_{v_3}-q_0q_{v_1}) \\
	2(q_{v_1}q_{v_3}-q_0q_{v_2}) & 2(q_{v_2}q_{v_3}+q_0q_{v_1}) & q_0^2-q_{v_1}^2 -q_{v_2}^2+q_{v_3}^2
	\end{bmatrix}.
	\end{equation*}
	In what follows we present how to get the swimmer frame dynamics in the laboratory frame. Coordinates $X$ in the laboratory frame are expressed in the swimmer frame by $[0,x] = q*[0,X]*\bar{q}$. Performing time differentiation we get
	\begin{equation*}
	\frac{\mathrm{d} [0,x]}{\mathrm{d} t} = \frac{\mathrm{d} q}{\mathrm{d}t}*[0,X]*\bar{q} + q*[0,X]*\frac{\mathrm{d}\bar{q}}{\mathrm{d}t} = \frac{\mathrm{d} q}{\mathrm{d}t}*\bar{q}*[0,x] + [0,x]*q*\frac{\mathrm{d}\bar{q}}{\mathrm{d}t}.
	\end{equation*}
	Since
	\begin{equation*}
		\frac{\mathrm{d}[0,x]}{\mathrm{d}t} = [0,2v \wedge x],
	\end{equation*}
	in our case $2[0,v] = 2\frac{\mathrm{d}q}{\mathrm{d}t}*\bar{q} = [0,R(q)\omega]$. Then, multiplying both sides of the equality by $q$ gives the following dynamics for the quaternion
		\begin{equation}
		\frac{\mathrm{d}q}{\mathrm{d}t} = \frac{1}{2}[0,R(q)\omega]*q.
		\label{Equation:quaternion}
		\end{equation}
	\section{Implementation}
	\label{Section:Implementation}
	The implementation part combined the FEEL++ finite element library \cite{FEEL2013} to solve the Stokes problems with an ODE solver for the solution of the rigid body motion. Matlab and Paraview were used for the post-processing. 
	
	Since the methodology giving the rigid body problem was posed in the swimmer's reference frame and the swimmer's motion in the computational domain produced strong deformations on the underlying triangular mesh, at every time-step the domain was displaced and remeshed. This procedure was case specific, since both the shape and the deformation velocity on the swimmer's boundary were externally imposed.
	
	Thanks to FEEL++, a Stokes solver with a syntax closely resembling the finite element variational formulation is achieved. Moreover, implementation of the finite elements was strongly simplified: Taylor-Hood spaces of order $P^2/P^1$ are chosen, but also spaces of higher order (ex. $P^3/P^2$) could be supported easily. 
	
	The solution of the ODEs \eqref{RigidMotion} describing rigid-body motion was embedded in the code through an external ODE solver in C++. A Runge-Kutta scheme was performed. At the end of each time-step, before the new rotation matrix is computed, the quaternion needs to be normalized to have unit norm: the propagation of numerical errors does not preserve the  norm of the quaternion. Remark that, if this renormalization is not performed, the rotation matrix won't be orthogonal and the body would suffer a deformation.
	\section{Numerical tests}
	\label{Section:Experiments}
	\subsection{Validation of fluid solver}
\begin{figure}[h]
		\begin{subfigure}{0.45\textwidth}
			\includegraphics[width=\textwidth]{\FiguresPath/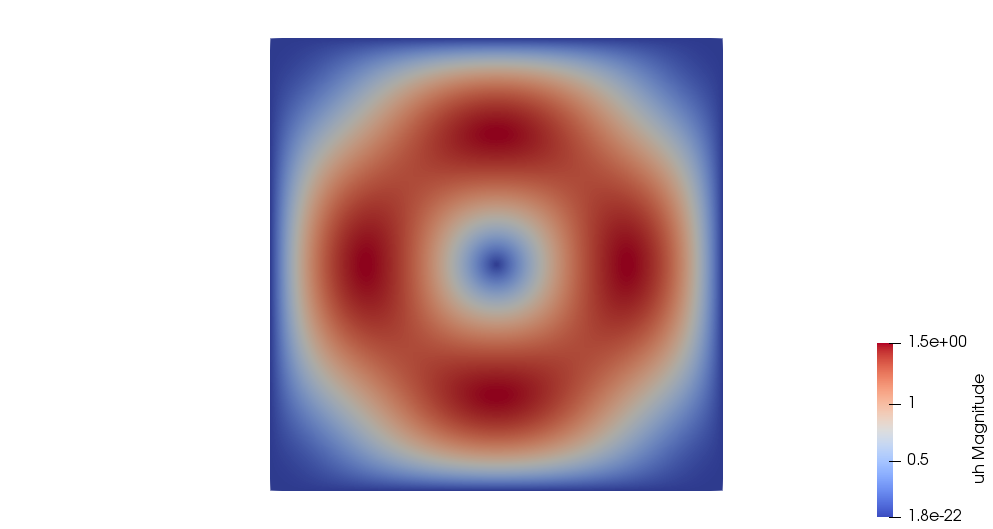}
			\caption*{Velocity magnitude}
		\end{subfigure}	
		\begin{subfigure}{0.45\textwidth}
			\includegraphics[width=\textwidth]{\FiguresPath/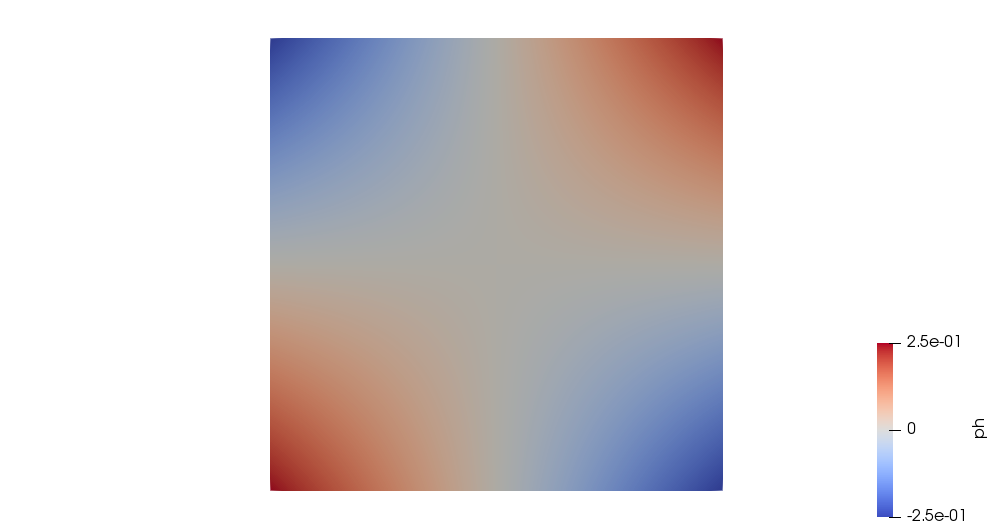}
			\caption*{Pressure}
		\end{subfigure}	
		\caption{Analytical solution of the Stokes problem \cite{Bercovier1979}.}
		\label{Fig:Bercovier-Engelman}
	\end{figure}

\begin{figure}[h]
		\begin{subfigure}{0.3\textwidth}
			\includegraphics[width=\textwidth]{\FiguresPath/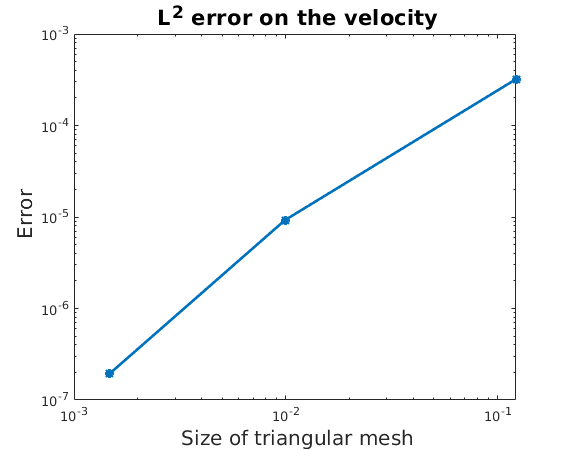}
		
		\end{subfigure}	
		\begin{subfigure}{0.3\textwidth}
			\includegraphics[width=\textwidth]{\FiguresPath/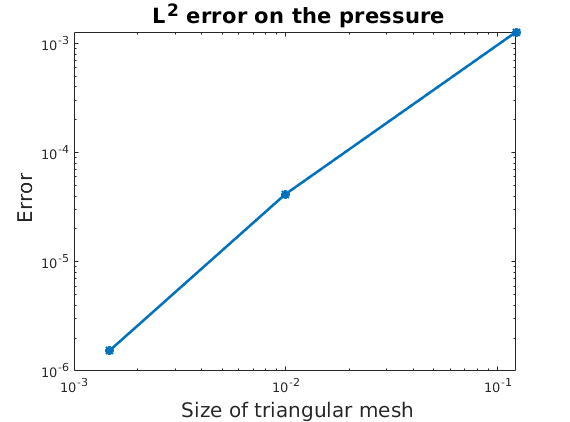}
		\end{subfigure}	
		\begin{subfigure}{0.3\textwidth}
			\includegraphics[width=\textwidth]{\FiguresPath/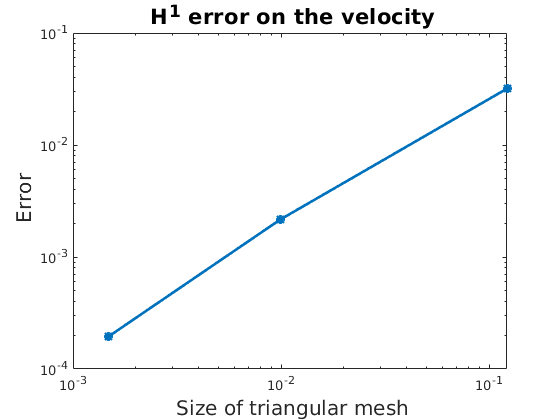}
		\end{subfigure}
		\caption{Error convergence of the numerical  solution given by our solver. The errors are in good agreement with the classical errors given by the finite element method: the experimental order of convergence for the velocity in $\mathrm{L}^2$ norm is approximately 3, and in $\mathrm{H}^1$ norm it's approximately 2. The experimental order of convergence for the $\mathrm{L}^2$ norm of the pressure is approximately 2.}
		\label{figure:Stokes_convergence}
	\end{figure}

	To validate the fluid solver presented in this paper, we compared the analytical solution of Stokes problem presented in \cite{Bercovier1979} to the numerical results given by our simulations. The analytical solution to Stokes problem on the unit square
	\begin{equation}
	\left\{
	\begin{aligned}
	&-\Delta u +\nabla p  = f \quad &\textrm{on $[0,1]\times[0,1]$}, \\
	&\nabla \cdot u = 0 \quad &\textrm{on $[0,1]\times[0,1]$},
	\end{aligned}
	\right.
	\label{eq:Stokes_on_square}
	\end{equation}
with forcing function
	\begin{equation*}
	f(t,(x,y)) = 
	\begin{bmatrix}
	128[x^2(x-1)^2 12(2y-1) + 2(y-1)(2y-1)y(12x^2-12x+2)] + y-0.5\\
	-128[y^2(y-1)^2 12(2x-1) + 2(x-1)(2x-1)x(12y^2-12y+2)]+x-0.5
	\end{bmatrix},
	\end{equation*}
	Dirichlet boundary conditions derived from the exact solution presented below and zero-mean pressure is given by
	\begin{equation}
		u(t,(x,y)) = 
		\begin{bmatrix}
			-256y(y-1)(2y-1)x^2(x-1)^2\\
			256x(x-1)(2x-1)y^2(y-1)^2
		\end{bmatrix},
		\label{eq:u_BE}
		\end{equation}
		\begin{equation*}
		p(t,(x,y)) = (x-0.5)(y-0.5).	\end{equation*}
	Figure \ref{Fig:Bercovier-Engelman} presents this latter Stokes solution. We solve the Stokes equations using the pair $(X_h^2,M_h^1)$ for the mixed finite elements.
	Figure \ref{figure:Stokes_convergence} shows a third order convergence of the $\mathrm{L}^2$ norm of the velocity and a second order convergence for both the $\mathrm{L}^2$ norm of the pressure and the $\mathrm{H}^1$ norm of the velocity are obtained. These results are those we expected from the finite element method theory, see for instance \cite[Thm 4.3, pag 181]{Girault1986}.
	
	Remark that, even though the solver is made for the 2D case, Feel++ allows to easily obtain a three dimensional Stokes solution when providing an appropriate 3D mesh.
	\subsection{Validation of rigid body solver}
	In order to validate the rigid body solver, we considered the results presented in \cite{Languski1993}. The solution to
	\begin{subequations}
	\begin{align}
	& \frac{\mathrm{d}X}{\mathrm{d}t} = R(q)V, \\
	& \frac{\mathrm{d}q}{\mathrm{d}t} = \frac{1}{2} [0,R(q)\omega]*q, \label{eq:subeq2} \\
	& \frac{\mathrm{d}V}{\mathrm{d}t} = F/m,\\
	& I\frac{\mathrm{d}\omega}{\mathrm{d}t}+ \omega \wedge I\omega = T,\label{eq:subeq4}
	\end{align}
	\label{RigidMotion_validation} 		\end{subequations}
	complemented with adequate initial conditions, is searched. 
	
	Let us remark that the exactness of the translational part of the rigid body motion is easy to check. As a result, we will consider only equations \eqref{eq:subeq2} and \eqref{eq:subeq4}.
	
\begin{figure}[h]
		\centering
		\includegraphics[width=0.47\textwidth]{\FiguresPath/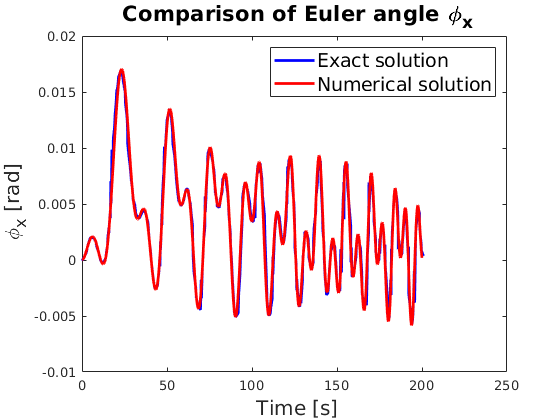}
		\includegraphics[width=0.45\textwidth]{\FiguresPath/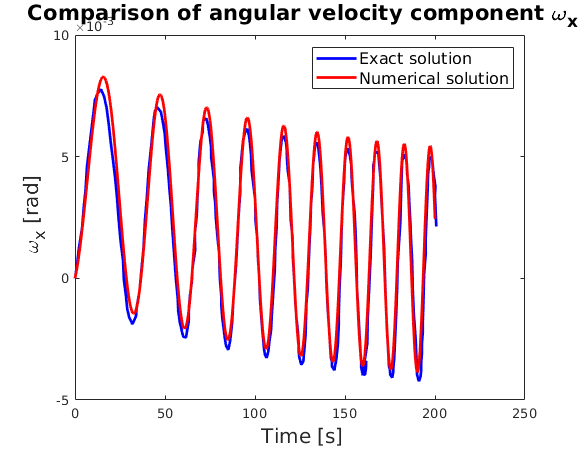}		
		\caption{Comparison of the sampled solution from \cite{Languski1993} and our numerical solution from equations \eqref{RigidMotion_validation}. On the left, the results on Euler angle $\phi_x$ are compared. On the right, the results on the component $\omega_x$ of the angular velocity are compared.}
		\label{figure:rigid_body}
	\end{figure}

	We validate the rotational part of the rigid body motion by extracting a set of sample points from the graphics of \cite{Languski1993} and compared qualitatively such results with ours. The case we will consider is the following: the time-varying torque in the body frame will be given by
	\begin{equation*}
		T = 
		\begin{bmatrix}
		1.0+2.7	\times 10^{-2}t-2.4\times 10^{-4}t^2 + 5.7\times 10^{-7}t^3 \\
		-1.5-9.0\times 10^{-3} t + 1.2\times 10^{-4}t^2-3.0\times 10^{-7}t^3\\
		13.5
		\end{bmatrix}N\cdot m .
	\end{equation*}
	The body reference frame at $t=0$ coincides with the laboratory reference frame, which in terms of quaternions is expressed as $q(0)=(1,0,0,0)$. The components of the angular velocity $\omega$ at $t=0$ along the three principal axes are $(0,0,0.329)$$\,rad/s$, while the inertia matrix is given by
	\begin{equation*}
	I = 
		\begin{bmatrix}
		2985 & 0& 0 \\
		0 & 2729 & 0 \\
		0 & 0 & 4183
		\end{bmatrix} kg\cdot m^2
	\end{equation*}
	Since in \cite{Languski1993} the results about the frame orientation are given in terms of Euler angles $(\phi_x,\phi_y,\phi_z)$, we use the following relationship to obtain them from a quaternion $q=(q_0,q_1,q_2,q_3)$:
	\begin{equation*}
		\begin{bmatrix}
		\phi_x \\ \phi_y \\ \phi_z
		\end{bmatrix} = 
		\begin{bmatrix}
		\arctan \frac{2(q_0q_1 + q_2q_3)}{1-2(q_1^2+q_2^2)} \\
		\arcsin (2(q_0q_2-q_3q_1)) \\
		\arctan \frac{2(q_0q_3 + q_1q_2)}{1-2(q_2^2+q_3^2)}
		\end{bmatrix}
	\end{equation*}
	Figure \ref{figure:rigid_body} shows that our results are in agreement with the numerical results in \cite{Languski1993}.
		
	We remark that, even though we will be focusing on a two-dimensional simulation for the swimmer, the previous validation confirms that such solver can be used also on a 3D experiment. The resulting angular velocity is the one required in \eqref{RigidMotion} to obtain the body orientation.
	\subsection{Validation of force and torque computations}
	Given $(u,p)$ solution of Stokes equations, this part addresses the numerical computation of fluid forces and torques. Since the swimmer motion depends on them, it is fundamental to validate the techniques we use to compute them. Two methods are compared: the first method consists in approximating the integrals describing surface forces and torques
	\[
	\int_{\partial \Omega} \sigma(t,x) \vec{n}(t,x) \, \mathrm{d}S(x), \qquad
	\int_{\partial \Omega} \sigma(t,x)\vec{n}(t,x)\wedge x \, \mathrm{d}S(x),
	\] 
	via surface quadratures, while the second method consists in exploiting the variational formulation to resort to volume integrals 
	\begin{equation}
	\begin{aligned}
	&\int_{\partial \Omega} \sigma(t,x) \vec{n}(t,x) \, \mathrm{d}S(x) = \int_{\Omega} f \, \mathrm{d}x + \int_{\Omega} \sigma(t,x) : \nabla_x u  \, \mathrm{d}x,
	\\
	&\int_{\partial \Omega} \sigma(t,x)\vec{n}(t,x)\wedge x \, \mathrm{d}S(x) = \int_{\Omega} f\wedge x \, \mathrm{d}x + \int_{\Omega} (\sigma(t,x) : \nabla_x u) \wedge x \, \mathrm{d}x,
	\end{aligned}
	\label{Eq:Volume Integral}
	\end{equation}
	that are subsequently approximated via 3D volumic quadratures. Feel++ library was used to perform these computations. The numerical results will be compared to approximate and exact solutions. We considered two test cases: the first one is the 3D Stokes example of a translating and rotating sphere, for which the approximate solutions are provided in \cite[p.67-68]{happel1983low}, while the second one is a 3D example taken from \cite{Ethier1994}. In the moving sphere case, the experimental setting is composed of a sphere of radius $r=0.01$ placed in the middle of a cubic box of size $1$, $\mu=1$ and $U=\omega=2$. The approximated forces and torques, given by the formulas $F = -6\pi \mu r U$ and $T=-8\pi \mu r^3 \omega$, give a qualitative target for the numerically derived forces and torques.
	
	The results regarding the 3D Stokes case are collected in Tables \ref{Table:3DStokes-force} to \ref{Table:3DStokes-torque-p3p2}. Tables \ref{Table:3DStokes-force} and \ref{Table:3DStokes-torque} show the numerical results that were obtained when using $X_h^2/M_h^1$ finite elements on a piecewise linear approximation of the geometry: Table \ref{Table:3DStokes-force} shows that the surface quadratures give better approximation of forces while Table \ref{Table:3DStokes-torque} show that torques are better approximated when using volumic quadratures. Tables \ref{Table:3DStokes-force-p3p2} and \ref{Table:3DStokes-torque-p3p2} show the numerical results that were obtained when $X_h^3/M_h^2$ finite elements on a piecewise quadratic approximation of the geometry were used: in this case volumic quadrature performs better than surface quadratures in both cases.
	\begin{center}
		\begin{table}[h]
		\begin{tabular}{||c c c c || c|} 
			\hline
			$h_{min}$ & $h_{max}$ & Surface quad. & Volume quad. & Stokes \\ [0.5ex] 
			\hline \hline
			0.00099 & 0.1958 & -0.397501 & -0.395855 & -0.376991 
			\\
			\hline
			0.00052 & 0.1546 & -0.394835 & -0.395245 & -0.376991
			\\
			\hline
			0.00036 & 0.1211 & -0.393037 & -0.393478 & -0.376991 
			\\
			\hline
		\end{tabular}
		\caption{Forces comparison for the translating and rotating sphere with piecewise linear approximation of the geometry and $X_h^2/M_h^1$ finite elements. In the last column, the numerical result of the Stokes formula $F = -6\pi \mu r U$ is reported.}
		\label{Table:3DStokes-force}
		\end{table}
		\begin{table}[h]
			\begin{tabular}{||c c c c || c|} 
				\hline
				$h_{min}$ & $h_{max}$ & Surface quad. & Volume quad. & Stokes \\ [0.5ex] 
				\hline\hline
				
				0.00099 & 0.1958 & $-4.97575\cdot 10^{-5}$ & $-5.0582\cdot 10^{-5}$ & $-5.02655 \cdot 10^{-5}$\\
				
				\hline
				0.00052 & 0.1546 & $-4.96165\cdot 10^{-5}$ & $-4.98904\cdot 10^{-5}$ & $-5.02655 \cdot 10^{-5}$
				
				\\
				\hline
				0.00036 & 0.1211 & $-5.03447\cdot 10^{-5}$ & $-5.03904 \cdot 10^{-5}$ &	$-5.05423 \cdot 10^{-5}$
				
				\\
				\hline
			\end{tabular}
			\caption{Torques comparison for the translating and rotating sphere with piecewise linear approximation of the geometry and $X_h^2/M_h^1$ finite elements. In the last column, the numerical result of the Stokes formula  $T=-8\pi \mu r^3 \omega$ is reported.}
			\label{Table:3DStokes-torque}
		\end{table}	
	\end{center}
	\begin{center}
		\begin{table}[h]
			\begin{tabular}{||c c c c|| c|} 
				\hline
				$h_{min}$ & $h_{max}$ & Surface quad. & Volume quad. & Stokes \\ [0.5ex] 
				\hline\hline

				\hline
				0.00099 & 0.1958 & -0.39264 & -0.392629 & -0.376991 
				
				\\ 
				\hline
				0.00052 & 0.1546 & -0.392439 & -0.392408 &-0.376991 
				
				\\
				\hline
			\end{tabular}
			\caption{Forces comparison for the translating and rotating sphere with piecewise quadratic approximation of the geometry and $X_h^3/M_h^2$ finite elements. In the last column, the numerical result of the Stokes formula $F = -6\pi \mu r U$ is reported.}
			\label{Table:3DStokes-force-p3p2}
		\end{table}
		\begin{table}[h]
			\begin{tabular}{||c c c c || c|} 
				\hline
				$h_{min}$ & $h_{max}$ & Surface quad. & Volume quad. & Stokes \\ [0.5ex] 
				\hline\hline
				0.00099 & 0.1958 & $-5.02279\cdot 10^{-5}$ & $-5.02783\cdot 10^{-5}$ & $-5.02655 \cdot 10^{-5}$\\
				
				\hline
				0.00052 & 0.1546 & $-5.02507\cdot 10^{-5}$ & $-5.02577\cdot 10^{-5}$ & $-5.02655 \cdot 10^{-5}$
				
					\\ 
			\hline
			\end{tabular}
			\caption{Torques comparison for the translating and rotating sphere with piecewise quadratic approximation of the geometry and $X_h^3/M_h^2$ finite elements. In the last column, the numerical result of the Stokes formula  $T=-8\pi \mu r^3 \omega$ is reported.}
			\label{Table:3DStokes-torque-p3p2}
		\end{table}	
	\end{center}
	The second benchmark that was considered regarded the exact computation of forces from an exact and fully three-dimensional solution to Navier-Stokes equations (where we neglected the time variation)  \cite{Ethier1994} 
	\begin{equation}
		u = -a
		\begin{bmatrix}
		e^{ax}\sin(ay+dz)+e^{az}\cos(ax+dy)\\
		e^{ay}\sin(az+dx)+e^{ax}\cos(ay+dz)\\
		e^{az}\sin(ax+dy)+e^{ay}\cos(az+dx)
		\end{bmatrix},
		\end{equation}
		\begin{multline}
		p = -\frac{a^2}{2}[e^{2ax}+e^{2ay}+e^{2az}+2\sin(ax+dz)\cos(az+dx)e^{ay+az} \\+2\sin(ay+dz)\cos(ax+dy)e^{az+ax}+2\sin(az+dx)\cos(ay+dz)e^{ax+ay}]
		\end{multline}
		in the cube $[-1,1]^3$, with $a=\frac{\pi}{4}$  and $d=\frac{\pi}{2}$.
		
\begin{figure}
		\begin{subfigure}{0.3\textwidth}
			
			\includegraphics[width=\textwidth]{\FiguresPath/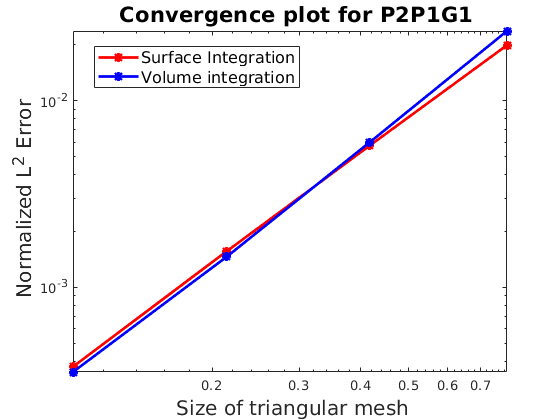}
			
		\end{subfigure}	
		\begin{subfigure}{0.3\textwidth}
			
			\includegraphics[width=\textwidth]{\FiguresPath/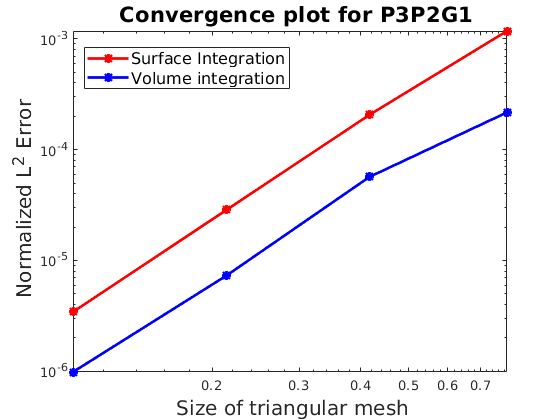}
		\end{subfigure}	
		\begin{subfigure}{0.3\textwidth}
			
			\includegraphics[width=\textwidth]{\FiguresPath/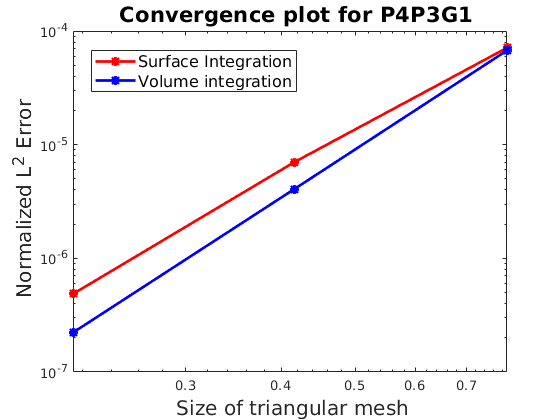}
			\end{subfigure}
			\caption{Comparison in $\mathrm{L}^2$ norm between exact and approximated forces, normalized by the $\mathrm{L}^2$ norm of the exact force. These pictures show that volume quadratures give a better approximation of forces and torques as the mesh size decreases. All the figures were obtained using a linear approximation of the geometry (G1). The left figure shows the results for $X_h^2/M_h^1$ finite elements. The center figure shows the results for $X_h^3/M_h^2$ finite elements. The right figure shows the results for  $X_h^4/M_h^3$ finite elements. }
			\label{Fig:Quadrature}			
	\end{figure}

	The exact forces were calculated using the symbolic mathematics package \emph{Sympy}, and they were compared to forces issued from surface and volume quadratures. The results in Figure \ref{Fig:Quadrature} show that, after fixing a high quadrature order for both surface and volume quadratures (13 in our case), forces are better approximated when using volume quadratures on \eqref{Eq:Volume Integral} than surface quadratures.
	These two benchmarks encouraged us to use the volume formulation to compute forces and torques for the Scallop simulation in the next section.

	\subsection{Scallop theorem}
\begin{figure}
		\begin{subfigure}{0.2\textwidth}
			\includegraphics[width=\textwidth]{\FiguresPath/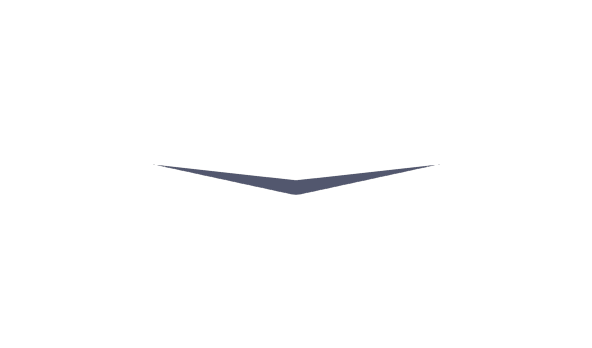}
			
		\end{subfigure}	
		\begin{subfigure}{0.2\textwidth}
			\includegraphics[width=\textwidth]{\FiguresPath/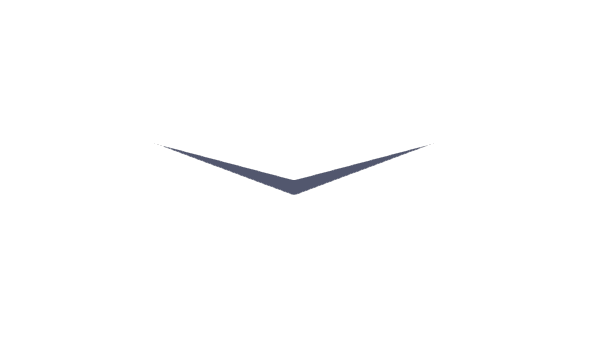}
		\end{subfigure}	
		\begin{subfigure}{0.2\textwidth}
			\includegraphics[width=\textwidth]{\FiguresPath/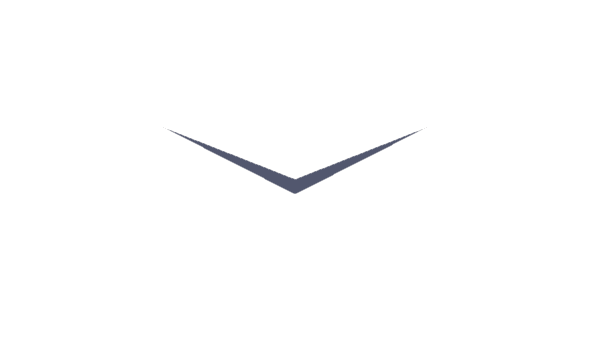}
		\end{subfigure}
		\begin{subfigure}{0.2\textwidth}
			\includegraphics[width=\textwidth]{\FiguresPath/Luca1.png}
		\end{subfigure}
		\caption{Screenshots of the scallop's body configurations over a period. The opening and closure of the swimmer's valves produce the motion. Since it has only one degree of freedom, the scallop does not experience a net advancement over a period.}
		\label{Fig:ScallopOpening}			
	\end{figure}

	In this section we will address the simulation of the so called ``Scallop theorem". This theorem was first presented by Purcell in \cite{Purcell1977}, and later generalized in \cite{Chambrion2010} for other types of swimmers, highlighting the underlying geometrical properties of the problem. Swimming is characterized by periodical strokes, and the scallop theorem specifies when they do not provide a net advancement. Scallop theorem, as in \cite{Purcell1977}, was specialized to Stokes flows, and stated that a swimmer which performs time-reversible shape changes will not be able to advance over a period. A particular case, from which the theorem takes his name, is the case of one-hinged swimmers. Having just one degree of freedom, swimmers like the scallop are forced to perform time-reversible shape changes and are not able to obtain a net advancement over a period of time. Figure \ref{Fig:ScallopOpening} shows an example of these shape changes.
	
	The aim of this part is to verify that our numerical solver is able to reproduce the scallop theorem. Problem \eqref{eq:Stokes} is solved and coupled with \eqref{RigidMotion} to obtain the swimmer's displacement and change in orientation. Here $\Omega^0$ is the complement in $[0,0.1]\times[0,0.1]$ of the scallop and the deformation velocity on the swimmer's boundary point $x$ is given by $v(t,x)=\frac{\mathrm{d}A}{\mathrm{d}t}(t)x$, where $A(t)$ is the rotation matrix imposing the desired opening and closing of the valves.
		
	We imposed on the boundary of the swimmer a displacement field which caused its hinges to open and close in a periodical fashion. We also allow a slight deformation of the body, to see whether any changes can be experienced with respect to the invariance predicted in \cite[Section 3.5]{Chambrion2010}. 
	
	The computational domain coincides with $\Omega^0$, hence it consists in a pierced square whose hole has the shape of the swimmer. The mesh is conformal to the swimmer's geometry and the approximation is linear: interior and boundary edges are straight.
	
\begin{figure}
		\begin{subfigure}{0.45\textwidth}
		
			\includegraphics[width=\textwidth]{\FiguresPath/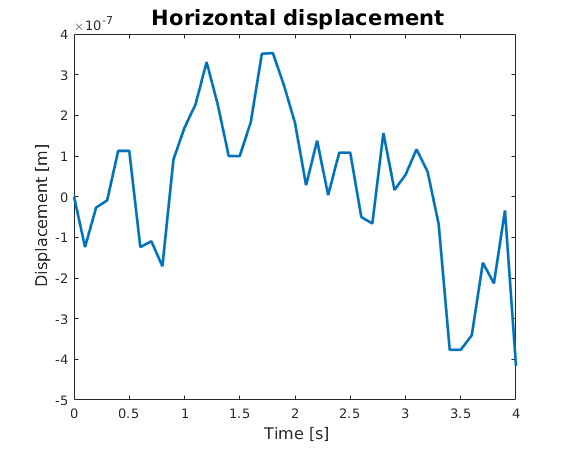}			
		\end{subfigure}	
		\begin{subfigure}{0.45\textwidth}
			
			\includegraphics[width=\textwidth]{\FiguresPath/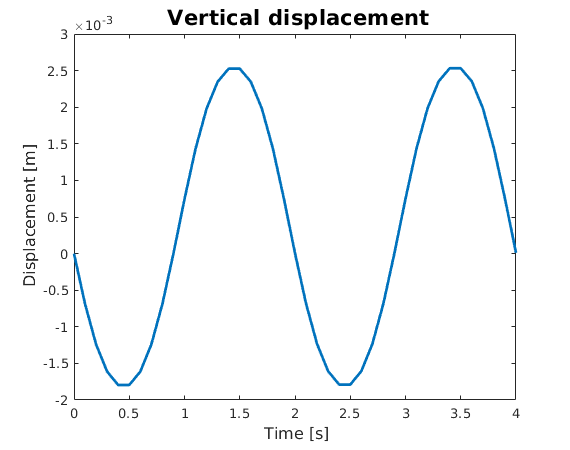}			
		\end{subfigure}	
		\caption{Position of the swimmer's center of mass in the laboratory reference frame. In this case the swimmer is located far from the no slip boundaries of the computational domain. The figure on the left shows its horizontal displacement, along the $x$ axis: it's globally negligible and imputable to numerical errors. The figure on the right shows the swimmer's vertical displacement: it is periodic of period $T=2$, and the net motion value on a period is also imputable to numerical errors.}
		\label{figure:scallop}
	\end{figure}

	Figure \ref{figure:scallop} shows that the scallop theorem is not exactly satisfied due to the domain discretization, but the errors over a period, and the errors cumulated between two consecutive periods, are negligible. We chose a period $T = 2\, s$ and a swimmer with total valve opening equal to $14 \, m m$. The smallest triangulation edge size in is of $0.025 \, m m$. 

	
\begin{figure}
			\begin{subfigure}{0.45\textwidth}
			
				\includegraphics[width=\textwidth]{\FiguresPath/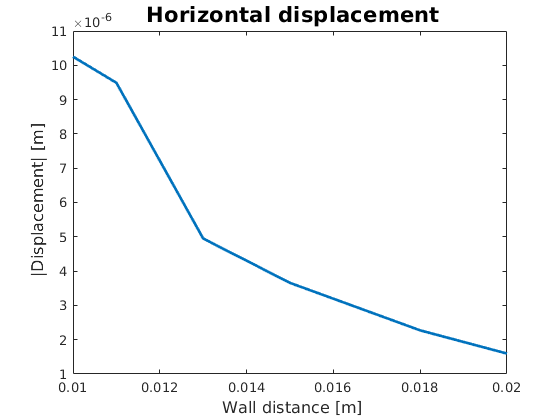}			
			\end{subfigure}	
			\begin{subfigure}{0.45\textwidth}
				
				\includegraphics[width=\textwidth]{\FiguresPath/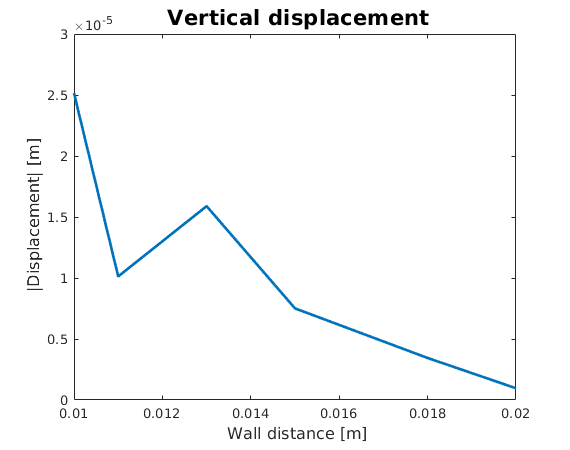}			
			\end{subfigure}	
			\caption{Absolute value of the net displacement of the swimmer's center of mass, as a function of its distance to a wall. When the swimmer departs closer to the no slip boundary, it experiences a larger displacement over a time period, both in the horizontal and vertical direction.}
			\label{figure:scallop_wall}
		\end{figure}

	We also located the scallop close to a no slip boundary of the computational domain. In that case the scallop theorem is not satsfied since the wall breaks the symmetry of the whole system. Figure \ref{figure:scallop_wall} shows the displacement over a complete stroke period as a function of the proximity to the wall. In that case we noticed that the displacement at the end of the period was sensibly larger than what Figure \ref{figure:scallop} showed, and it exemplifies both the attractive effect that walls have on swimmers (Figure \ref{figure:scallop_wall} right) and the enhanced displacement parallel to no slip boundaries (Figure \ref{figure:scallop_wall} left). Let us quote the recent study \cite{Bagagiolo2017} that overcomes the scallop theorem by opening and closing the shell with different speeds.
	
	The pseudocode for the simulation of the swimmer motion is presented in Algorithm \ref{Algo:Code}.

\begin{algorithm}
	\caption{Solution procedure for the swimming problem}
	\label{Algo:Code}
	\begin{algorithmic}
		\STATE Choose $\Delta t$, $d$, $T_{final}$
		\STATE Initialize $X, q, V, \omega$
		\IF{$d=2$}
		\STATE $N_{Stokes} =4$
		\ELSIF{$d=3$}
		\STATE $N_{Stokes}=7$
		\ENDIF
		\WHILE{$time \le T_{final}$}
		\STATE Displace the swimmer according to $X, q, V, \omega$ and Remesh
		\FOR{$1 \le i \le N_{Stokes}$}
		\STATE Solve Stokes problem $i$ from subsection \ref{Subsection:RigidBody}
		\STATE Compute surface forces and torques issued from Stokes problem $i$
		\ENDFOR
		\STATE Find $V,\omega$ by using \eqref{Eq:M} and \eqref{Eq:N}
		\STATE Find $X,q$ by solving \eqref{RigidMotion}
		\STATE $time \leftarrow$ $time + \Delta t$
		\ENDWHILE
	\end{algorithmic}
\end{algorithm}
	\section{Conclusion and perspectives}
	In this note we present a solver able to provide the dynamics of a derforming swimmer moving into a Stokes flow. Here all the framework could fit a 3D model but the numerical results are given in 2D. All the numerical computations of this note relied on the finite element library Feel++. We applied a method proposed in \cite{Munnier2012} to extract rigid body motion from computation of fluid stresses. We validated our code and model using different benchmarks. The goal of this project was to lay the foundations of a numerical method which allowed the simulation of swimmers with more complex geometries.
	
	The passage to several swimmers or to 3D are beyond the scope of this paper. The technique could be applied to several swimmers, but the basis onto which the Stokes solution has to be expanded is larger and the computational cost might increase noticeably. Another issue is that, in the Stokes regime, forces can become infinite when two swimmers are close, and collision avoidance strategies would need to be implemented.  
	\bibliography{swimmer.bib}
	\bibliographystyle{siam}
\end{document}